\documentclass{amsart}
\usepackage{amssymb}
\usepackage{color}
\usepackage{graphicx}
\usepackage{amsmath,amscd}
\usepackage{amsthm}
\usepackage{float}
\usepackage{bm}
\usepackage{tikz}
\usepackage{tikz-cd}
\usepackage{enumitem}
\usepackage{epstopdf}
\usepackage{stackengine}
\usepackage{tikz-cd}
\usepackage{leftindex}
\usepackage{mathtools}
\usepackage[export]{adjustbox}
\usetikzlibrary{arrows.meta}
\usepackage[colorlinks=true,
            linkcolor=blue!80!black,
            citecolor=red!80!white,
            urlcolor=cyan!60!black]{hyperref}
\usepackage[pagewise]{lineno}

\newtheorem{theorem}{Theorem}[section]
\newtheorem{proposition}[theorem]{Proposition}
\newtheorem{definition}[theorem]{Definition}
\newtheorem{lemma}[theorem]{Lemma}

\newtheorem{remark}[theorem]{Remark}
\newtheorem{question}[theorem]{Question}

\renewcommand{\S}{\mathbb{S}}

\title{New Surgery Presentations and Applications to Open Book Embeddings}

\date{\today}



\thanks{Authors' names are listed alphabetically, reflecting equal contributions from each author to this research paper.}

\author{Sanjay Kute}
\address{Indian Institute of Technology Madras,
	IIT PO. Chennai- 600036, Tamil Nadu, India.ORCID ID: 0009-0005-4528-2763}
\email{sanjaykutemjt@gmail.com} 

\author{Dishant M. Pancholi}
\address{The Institute of Mathematical sciences,
	X6VW+FP7, 4th Cross St, CIT Campus, Tharamani, Chennai-600113 , Tamil Nadu India.}
\address{Homi Bhabha National Institute, Training School Complex, Anushaktinagar, Mumbai 400094, India}
\email{dishant@imsc.res.in} 

\author{Suhas Pandit}
\address{Indian Institute of Technology Madras,
	IIT PO. Chennai- 600036, Tamil Nadu, India. ORCID ID:0000-0001-7463-8523 }
\email{suhas@smail.iitm.ac.in}

\begin{document}

\begin{abstract} 
    We show that every closed, orientable $3$-manifold can be obtained from $\S^3$ by surgery along a link whose components are unknots and whose surgery coefficients are $0$ or $\pm\frac{1}{2}$. Using this result, we prove that every closed, orientable $3$-manifold $M$ admits an embedding into $\mathbb{S}^5 \subset \mathbb{C}^3$ such that the map
\[
    \pi \colon \mathbb{C}^3 \setminus \{z_3=0\} \longrightarrow \mathbb{S}^1, \qquad \pi(z_1,z_2,z_3)=\frac{z_3}{|z_3|},
\]
    induces an open book decomposition of $M$ with binding $M\cap\{z_3=0\}$. We conclude by discussing the construction of potentially new invariants of closed, orientable $3$-manifolds arising from open book decompositions and their embeddings.
\end{abstract}

\maketitle
\section{Introduction}\label{sec:intro}
    Dehn surgery is an important tool for the study of $3$-dimensional topology. Lickorish~ \cite{Lickorish1962} and Wallace~\cite{Wallace1960} proved that every closed, orientable $3$-manifold can be obtained by performing integer surgery along a link in $\mathbb{S}^3$. Kirby~\cite{Kirby1978} proved that any two framed links in $\mathbb{S}^3$ yielding orientation-preservingly diffeomorphic $3$-manifolds are related by a finite sequence of Kirby moves: (i) \emph{blow-up/blow-down}, adding or removing an unlinked unknot with framing $\pm 1$, and (ii) \emph{handle slides}, sliding one component over another.
    
    Surgery presentations of $3$-manifolds are central to  the study
    of geometric structures  on $3$--manifolds. For examples, the construction of hyperbolic structures~\cite{Thurston1980}, contact structures~\cite{Geiges2008},  open book decompositions~\cite{ROLF}.
    
    In this article, we give a new surgery description of closed, orientable $3$-manifolds, as stated in the following theorem.
\begin{theorem}\label{thm:binary_surgery_link}
    Every closed, oriented $3$-manifold can be obtained by surgery along a link $L$ in $\S^3$ which satisfies the following properties:
\begin{enumerate}
    \item each component of the link $L$ is an unknot,
    \item the linking number of any pair of components is $1$, $0$, or $-1$,
    \item the surgery coefficient of each component is $0$ or $\pm\frac{1}{2}$.
\end{enumerate}
\end{theorem}
    We give an interesting application of this surgery description to open book embeddings of $3$-manifolds into $\mathbb{S}^5$ in Theorem~\ref{thm:planar_page_main} and use it to describe potentially new invariants of $3$--manifolds.

    Recall that an open book decomposition $(B,\pi)$ of a closed $n$-manifold $M$ 
    is an expression  of the complement of a codimension $2$ submanifold $B$ of  $M$ as a fiber bundle $\pi\colon M\setminus B\to \mathbb{S}^1$. Fibers of $\pi$ -- or sometimes
    closures of fibers in $M$ -- are called \emph{pages}, while the submanifold $B$ is called the \emph{binding} of
    the open book $(B, \pi)$. Alexander~\cite{alexander1923} proved that every closed, orientable $3$-manifold admits an open book decomposition. Giroux ~\cite{Giroux2002} showed that there is a one to one correspondence between the isotopy classes of co-oriented contact structures on a closed, oriented, smooth $3$-manifold and the open book decompositions of $M$ up to positive stabilizations. 

    Several invariants of  $3$-manifolds based on open book decompositions  have been studied, including page genus, binding components, and monodromy invariants. In particular, Etnyre and Ozbagci ~\cite{EtnyreOzbagci2008} introduced three invariants of contact structures: \emph{support genus}, the minimum genus of a supporting page; \emph{binding number}, the minimum number of binding components among supporting open books of minimal genus; and \emph{support norm}, the minimum negative Euler characteristic of a supporting page. 

    A smooth manifold $M$ with a given open book decomposition $(B, \pi)$ admits an \emph{open book embedding} in a manifold $N$  with an open book decomposition $(B^\prime, \pi^\prime)$, provided there is an embedding of $M$ in $N$ such that –as a submanifold of $N$– the given open book decomposition $(B, \pi)$ of $M$ is obtained via the intersection of the open book decomposition $(B^\prime, \pi^\prime)$ of $N$ with $M$. If this is the case, we sometime say that $(B, \pi)$ is a pull-back of the open book $(B^\prime, \pi^\prime)$.

    Using the surgery description given in Theorem~\ref{thm:binary_surgery_link}, we establish the existence of open book embeddings of closed, orientable $3$-manifolds into $\mathbb{S}^5$ with the trivial open book. In particular, we prove the following:
\begin{theorem}\label{thm:planar_page_main}
    Let $M$ be a closed, orientable $3$-manifold. There exists an embedding $f$ of $M$ into $\S^5$ such that the trivial open book decomposition of $\S^5$ induces a planar open book decomposition of $M$.
\end{theorem}
    Here by the trivial open book of $\mathbb{S}^5$, we mean the open book decomposition of $\mathbb{S}^5$ having pages diffeomorphic to $4$-balls and monodromy map the identity.  Although every open book decomposition of closed, orientable $3$-manifold  can be realized as a pull-back~\cite{BanerjeeLawandeRanaSaha2026}, to the best of our knowledge, there is no systematic method that allows one to prescribe or  control the topology of the binding and pages of the corresponding open book decomposition of the target manifold. In this article, we try to address this issue.
 
    Motivated by Theorem~\ref{thm:planar_page_main}, we introduce a potentially new invariant of closed, orientable $3$-manifolds, called the \emph{embedding binding number}. It is defined as the minimum number of binding components among all planar open book decompositions that admit an open book embedding into the trivial open book of $\mathbb{S}^5$.
   
    Using an ambient stabilization process, we establish the following:
\begin{theorem}\label{thm:connected_binding_main}
    Let $M$ be a closed, orientable $3$-manifold. There exists an embedding $f$ of $M$ into $ \S^5$ such that the trivial open book decomposition of $\S^5$ induces an open book decomposition on $M$ with connected binding.
\end{theorem}
    Theorem~\ref{thm:connected_binding_main} allows us to introduce  the notion of \emph{embedding genus} $\mathcal{EG}(M)$ which measures the minimum genus of a page among all open book decompositions of $M$ with connected binding that admit an open book embedding into the trivial open book of $\mathbb{S}^5$. Comparing this invariant with the \emph{minimal open book genus} $\mathcal{MG}(M)$ allows us to investigate how the existence of an open book embedding constrains the topology of the pages. Although $\mathcal{MG}(M)=\mathcal{EG}(M)$ is known when $\mathcal{EG}(M)\leq 2$, whether this equality holds in general remains an interesting open question. 
\bigskip

\section {Preliminaries}\label{sec:prelim}

    In this section, we quickly review the notions of open book decompositions and abstract open book decompositions associated to manifolds. In this article, all the manifolds are connected, oriented, and the diffeomorphisms are orientation preserving.

    Recall by the \emph{mapping class group} $\mathcal{M}(\Sigma, \partial \Sigma)$ of a compact, connected, oriented surface $\Sigma$ with non-empty boundary we mean the group of isotopy classes of orientation preserving diffeomorphisms of $\Sigma$ which are the identity near the boundary of $\Sigma$. It is a well-known fact due to Lickorish~\cite{Lickorish1962}, \cite{Lickorish1964}  that the mapping class group of a compact, connected, oriented surface $\Sigma$ is generated by a finite number of Dehn twists along certain non-separating simple closed curves  on $\Sigma$. Recall by a \emph{Dehn twist} along a simple closed curve $\gamma$ on a surface $\Sigma$, we mean a diffeomorphism $\tau_\gamma: \Sigma\to \Sigma$ supported in a regular neighborhood of $\gamma$ in $\Sigma$ obtained by cutting $\Sigma$ along $\gamma$, twisting by an angle $2 \pi$ to the right or left and then re-gluing.

\medskip
\subsection{Open Book Decompositions}\mbox{}

    We recall the definition of open book decompositions of closed, oriented manifolds.
\begin{definition}[Open Book Decomposition]\label{def: OBD} 
    An open book decomposition of a closed, oriented manifold $M$ is a pair $(B,\pi)$, where 
\begin{enumerate}
    \item $B$ is a closed, co-dimension $2$ submanifold of $M$ with trivial normal bundle,
    \item $\pi: M \setminus B \to S^1$ is a fiber bundle,
    \item there is a  trivialization of a tubular neighborhood $\mathcal{N}(B)\cong  D^2\times B$ of $B$ in $M$ such that the map $\pi:  D^2\times B \setminus (\{0\}\times B) \to S^1 $ sends $(r, \theta, x)$ to $\theta,$ where $x\in B$ and $(r,\theta)$ are the polar coordinates on $D^2$, and 
    \item for each $\theta\in S^1$, the closure of $\pi^{-1}(\theta)$ is a compact, connected, oriented manifold $\Sigma$ with $\partial \Sigma= B$.
\end{enumerate}
    Here, $B$ is called the binding, and $\Sigma$ is called the page of the open book $(B, \pi)$. 
\end{definition}
    If a manifold $M$ is equipped with an open book $(B,\pi)$, then we denote this by $M=\mathcal{O}b(B,\pi)$. Note that the fibration $\pi: M\setminus \mathcal{N}(B) \to S^1$ is determined by the fiber $\Sigma=\pi^{-1}(\theta)$ and the monodromy $\phi$ of the fibration $\pi$. The monodromy $\phi:\Sigma\to \Sigma$ is obtained by restricting a fiber preserving flow on $M$ which fixes the binding pointwise. Moreover, $\phi$ can be chosen to be the identity map in a collar of $\partial \Sigma$. Therefore, $$M\setminus \mathcal{N}(B) \cong \mathcal{M}T(\Sigma, \phi)= \frac{\Sigma\times [0,1]}{(x,1 )\sim (\phi(x),0)}.$$  From this one can easily see that there is a diffeomorphism $F:  M_{(\Sigma, \phi)} \to M $, where $$ M_{(\Sigma, \phi)} = \mathcal{M}T(\Sigma, \phi) \bigcup_{Id}  \partial \Sigma \times D^2,$$ such that
\begin{enumerate}
    \item the map $\pi\circ F: \mathcal{M}T(\Sigma, \phi) \to S^1$ is given by $(\pi\circ F)(y,\theta)=e^{i\theta}$, 
    \item the map $$\pi\circ F: \partial\Sigma \times D^2 \setminus (\partial \Sigma \times \{0\} ) \to S^1$$ is given by $(\pi\circ F)( x,r,\theta)=e^{i\theta},$ where $x\in \partial \Sigma $ and $(r,\theta)$ are the polar coordinates on $D^2.$
\end{enumerate}
    Note that the monodromy $\phi$ of the fibration is determined by the fibration up to isotopy and conjugation by a self-diffeomorphism of $\Sigma$. This motivates the notion of an abstract open book as defined below.
\begin{definition}[Abstract Open Book Decomposition]
    An abstract open book is a pair $(\Sigma, \phi)$ where
\begin{enumerate}
    \item $\Sigma$ is a compact, oriented surface with non-empty boundary, and
    \item $\phi : \Sigma \to \Sigma$ is a diffeomorphism such that $\phi$ is the identity in a collar neighborhood of $\partial\Sigma$.
\end{enumerate}
    The map $\phi$ is called the \emph{monodromy}.
\end{definition}
    Given an abstract open book $(\Sigma, \phi)$, we can associate a manifold $$M = \mathcal{M}T(\Sigma, \phi) \bigcup_{Id}  \partial \Sigma \times D^2.$$ In this case we denote the manifold $M$ by $M=\mathcal{A}ob(\Sigma,\phi).$ Now, one can easily see that given an open book of a closed manifold, we can associate an abstract open book and vice versa \cite{Etnyre2006}. Now, we have the following straightforward proposition.

\begin{proposition}\cite{Etnyre2006}
    Let $M=\mathcal{A}ob(\Sigma, \phi)$. Then $M$ is diffeomorphic to $M^\prime =\mathcal{A}ob(\Sigma, f \circ  \phi \circ f^{-1})$ for every diffeomorphism $f: \Sigma \rightarrow \Sigma.$
\end{proposition}
\medskip

\subsection{Stabilization of Open Books}\mbox{}

    Let $(B,\pi)$ be an open book decomposition of $M$ with page $\Sigma$ and associated abstract open book $(\Sigma,\phi)$. Let $C \subset \Sigma$ be a properly embedded simple arc. We say that a compact surface $\Sigma' \subset M$ is obtained from $\Sigma$ by plumbing a positive/negative Hopf band along $C$, if $\Sigma'=\Sigma\cup \mathcal{A}^{\pm},$ where $\mathcal{A}^{\pm}$ is an annulus embedded in $M$ satisfying the following properties:
\begin{enumerate}
    \item $\mathcal{A}^{\pm}\cap\Sigma$ is a tubular neighborhood of $C$ in $\Sigma$;
    \item the core curve of $\mathcal{A}^{\pm}$ bounds an embedded disk in $M\setminus\Sigma$, and the linking number of the two boundary components of $\mathcal{A}^{\pm}$ is equal to $\pm1$.
\end{enumerate}
\begin{figure}[ht]
    \centering
    \includegraphics[width=0.5\linewidth]{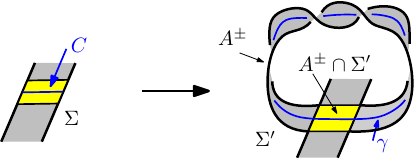}
    \caption{Plumbing of a positive Hopf band.}
    \label{Stabilization of open book}
\end{figure}
    According to a result of Stallings~\cite{JRS}, the surface $\Sigma'$ is also a page of an open book $(B',\pi')$ of $M$ whose associated abstract open book is $(\Sigma',\phi\circ\tau_{\gamma}^{\pm1}),$ where $\gamma$ is the core curve of the Hopf band $\mathcal{A}^{\pm}$.

    We say that the open book $(B',\pi')$ is obtained from $(B,\pi)$ by an $\mathcal{A}^{\pm}$-plumbing. A \emph{positive $(\text{negative})$ stabilization} of the open book $(B,\pi)$ is an open book $(B',\pi')$ that can be obtained from $(B,\pi)$ by a finite sequence of positive (negative) Hopf band plumbings.

\medskip
\subsection{Open Book Embeddings}\mbox{}

    Now, we recall the notion of open book embeddings of manifolds as open books.
\begin{definition}[Open Book Embedding]
    Let $M=\mathcal{O}b(B_1,\pi_1)$ and  $N=\mathcal{O}b(B_2, \pi_2)$ be two manifolds. We say that an embedding $f: M \hookrightarrow N$ is an open book embedding of $M=\mathcal{O}b(B_1,\pi_1)$  into $N=\mathcal{O}b(B_2, \pi_2)$ if $f$ embeds $B_1$ into $B_2$ and the following diagram commutes:
\begin{equation*}
\begin{CD}
    M \setminus B_1   @>f>>  N \setminus B_2\\
    @VV\pi_1V        @VV\pi_2V\\
    S^1     @>Id>>  S^1.
\end{CD}
\end{equation*}
\end{definition}
    Similarly, we can define an abstract open book embedding as follows:

\begin{definition}[Abstract Open Book Embedding]
    Let $M = \mathcal{A}ob \left(\Sigma_1, \phi_1\right)$ and $N = \mathcal{A}ob\left( \Sigma_2, \phi_2 \right)$ be two manifolds.  We say that there exists an abstract open book embedding of $M$ into $N$  provided there exists a proper embedding  $f$ of $\Sigma_1$ into $\Sigma_2$ such that $\phi_2  = f^{-1} \circ \phi_1 \circ f.$
\end{definition}
    It is clear from the definition that an abstract open book embedding produces an open book embedding for the associated open books and vice versa.
\bigskip

\section{A new surgery presentation of 3-manifolds}

\subsection{A Reivew of Rolfsen Twists}\mbox{}
    
    Rolfsen extended Kirby's calculus for integral framed links to links with rational framings. The key new move that Rolfsen introduced is a \emph{Rolfsen twist}~\cite{ROLF}.
\begin{figure}[ht]
    \centering
    \includegraphics[width=0.7\linewidth]{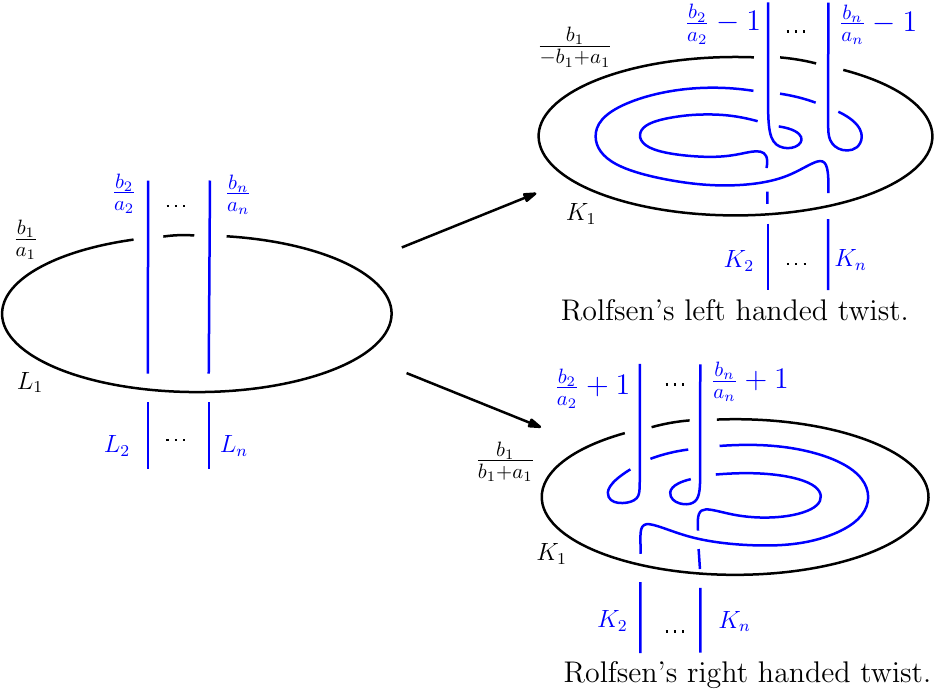}
    \caption{Rolfsen twists.}
    \label{Rolfsen_twist}
\end{figure}  

    In particular, Rolfsen proved that if $M$ is a closed, oriented $3$-manifold obtained by surgery along a framed link $L = L_1 \sqcup L_2 \sqcup  \cdots \sqcup L_m$ in $\S^3$, where the component $L_i$ has framing $\frac{b_i}{a_i}$ and $L_1$ is an unknot, then $M$ is obtained by surgery along a link  $\widehat{L}\subset \S^3$ with components $ L_1, K_2, \dots, K_m$, where for each $j = 2,\dots, m$, the component  $K_j$ is obtained from $L_j$ by twisting each strand of $L_j$ passing through a disk $D_1$ bounded by $L_1$ by $360^\circ$ rotation as depicted in Figure~\ref{Rolfsen_twist}. Furthermore, the new surgery coefficient of $L_1$ is $\frac{b_1}{a_1 \pm b_1}$ and for any $j \neq 1$ the surgery coefficient of $K_j$ is $(b_j/a_j) \pm \operatorname{lk}(L_1,L_j)^2$.
\medskip

\subsection{Surgery Description}\mbox{}   

    We are now ready to give a proof of Theorem~\ref{thm:binary_surgery_link}.
    We first prove the following lemma:
\begin{lemma}\label{lem:surgery_link}
    Let $M$ be a closed, oriented $3$-manifold obtained from $\S^3$ by surgery along a link $L = L_1 \sqcup \cdots \sqcup  L_m$. Suppose that the component $L_1$ is an unknot with an even integer surgery coefficient $n=2k$, where $k\in \mathbb{Z}$. Then $M$ can be obtained from $\S^3$ by surgery along a link $\widehat{L}$ which satisfies the following:
\begin{enumerate}
    \item $\widehat{L} = L \sqcup U_1\sqcup \cdots \sqcup U_{|k|}$, where each $U_i$ is an unknot,
    \item the surgery coefficient of $U_i$ is  $\pm\frac{1}{2}$,
    \item the surgery coefficient of $L_1$ is $0$,
    \item $\operatorname{lk}(U_i,L_1)=\pm1$ for each $i$, $\operatorname{lk}(U_i,L_j)=0$ for each $j\neq1$, and $\operatorname{lk}(U_i,U_l)=0$ for each $l\neq i$;, 
    \item the surgery coefficients of all remaining components $L_j$, $j\neq 1$, remain unchanged.
\end{enumerate}
\end{lemma}
\begin{proof}
    Assume that the surgery coefficient of $L_1$ is $n=2k$, where $k>0$. In this case, we first add an unknot $U_1$ with framing $-1$ inside a small product tubular neighborhood of $L_1$, as shown in Figure~\ref{Addition_of_unknot}. The unknot $U_1$ is chosen to be the meridian of $L_1$ with the opposite orientation, so that $\operatorname{lk}(U_1,L_1)=-1.$
\begin{figure}[ht]
    \centering
    \includegraphics[width=0.3\linewidth]{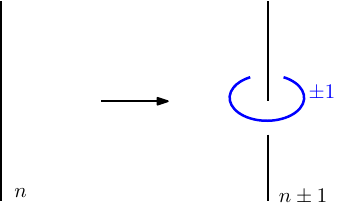}
    \caption{Addition of an unknot-blowing up.}
    \label{Addition_of_unknot}
\end{figure}

    By the Fundamental Theorem of Kirby calculus, the manifold $M$ can be obtained by surgery on the link $L\sqcup U_1.$ The surgery coefficient of $L_1$ changes from $2k$ to $2k-1$, while the surgery coefficients of the remaining components $L_j$, for $j\neq 1$, remain unchanged.

    Next, we perform a left-handed Rolfsen twist about the unknot $U_1$. By Rolfsen's theorem after the twist, the surgery coefficient of $L_1$ is $2k-2,$ and the surgery coefficient of $U_1$ is $-\frac{1}{2}.$ The surgery coefficients of all other components remain unchanged.

    We now repeat the same construction. Add another unknot $U_2$ with framing $-1$ as the oppositely oriented meridian of $L_1$, and perform a left-handed Rolfsen twist about $U_2$. The resulting surgery link is $L\sqcup U_1\sqcup U_2,$ for which the surgery coefficient of $L_1$ becomes $2k-4,$ while the surgery coefficient of $U_2$ becomes $-\frac{1}{2}.$ Again, the surgery coefficients of all remaining components are unchanged.

    Continuing this procedure inductively, after performing it $k$ times, the surgery coefficient of $L_1$ is $0$, the surgery coefficient of each unknot $U_i$ is $-\frac{1}{2}$ for $i=1,\dots,k$, and the surgery coefficients of all the other components of the link $L$ remain unchanged. 
    
    Similarly, if $k<0$, then by considering unknots $U_1,\dots, U_{|k|}$ as described above with framing $1$ and performing right-handed Rolfsen twists along each $U_i$, we get the surgery link $$\widehat{L} = L \sqcup U_1\sqcup \cdots \sqcup U_{|k|}$$ in $\S^3$ for the $3$-manifold $M$. The surgery coefficient of $L_1$ is $0$, the surgery coefficient of each unknot $U_i$ is $\frac{1}{2}$ for $i=1,\dots,{|k|}$, and the surgery coefficients of all the other components of the link $L$ remain unchanged. This establishes the lemma.
\end{proof}
\medskip

\begin{proof}[\textbf{A Proof of Theorem}~\ref{thm:binary_surgery_link}]
    Consider the surgery description of $M$ given by Birman and Powell~\cite{BirmanPowell1979}, as explicitly described by Hilden, Lozano, and Montesinos~\cite{HildenLozanoMontesinos1983}. They show that every closed, orientable $3$-manifold is obtained by surgery on a framed link  having each component unknotted with even integral framing. See Figure~\ref{The_surgery_link_given_by_Hilden_Lozano _and _Montesinos} for a pictorial description of one such link.
\begin{figure}[ht]
    \centering
    \includegraphics[width=0.8\linewidth]{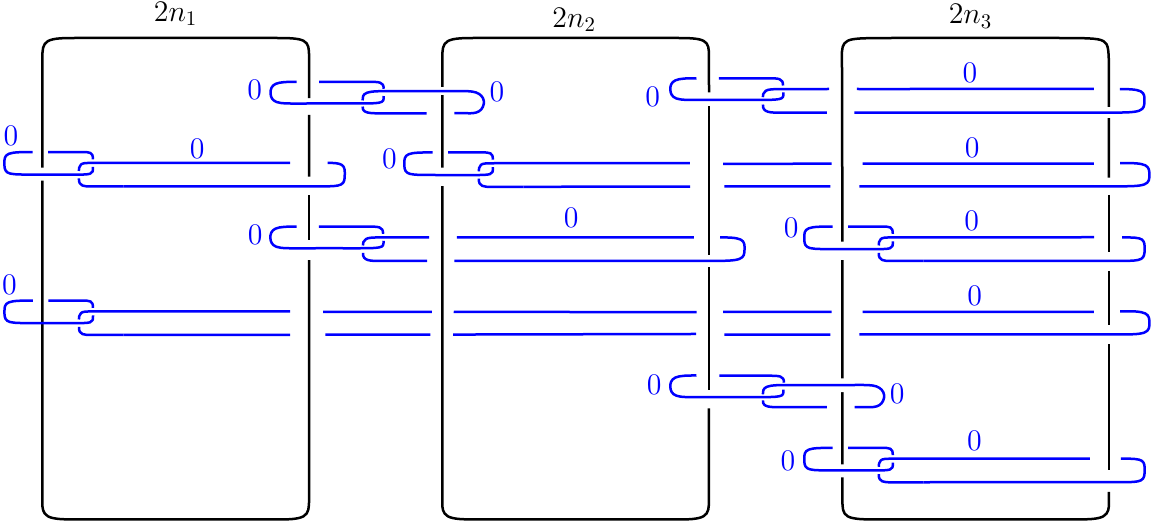}
    \caption{A typical surgery link $L$ with even integral framing coefficient as depicted by Hilden, Lozano, and Montesinos in
    \cite{HildenLozanoMontesinos1983}.}
    \label{The_surgery_link_given_by_Hilden_Lozano _and _Montesinos}
\end{figure}
    This implies that there exists a framed link $L$ in $\S^3$ such that $M$ is obtained from $S^3$ by surgery along the link $L$ and each
    component of $L$ is framed unknot with even framing. We apply Lemma~\ref{lem:surgery_link} individually to each component of $L$ whose surgery coefficient is non-zero. As a result, we obtain a surgery link $\widehat{L}$ in $\S^3$ whose surgery coefficients are all either $0$ or $\pm\frac{1}{2}$.
\end{proof}
    A framed link $L\subset \S^3$ satisfying the properties stated in Theorem~\ref{thm:binary_surgery_link} will be called a \emph{binary link} and the surgery coefficients of the components, the \emph{binary coefficients}.
\bigskip

\section{Embedded surfaces in \texorpdfstring{$\mathbb{S}^3$}{S3}}

    We need some basic facts about \emph{fiber surfaces} in $\S^3$. Let us begin with the definition of a fiber surface.
\begin{definition}[Fiber surface]
    A compact surface  $\Sigma$ embedded in $\S^3$ is said to be a fiber surface if there exists an open book of $\S^3$ with the binding $\partial \Sigma$ and $\Sigma$ is a page of this open book.
\end{definition}

    Clearly, by definition, the boundary of a fiber surface in $\S^3$ is a fibered link in $\S^3$. We recall the following result of Onaran from~\cite{Onaran2010}.
\begin{theorem}[\cite{Onaran2010}, Onaran]\label{thm:sinem}
    Any link in $\S^3$ lies on a fiber surface of a planar open book of $\S^3$.
\end{theorem}

    We recall the following result from ~\cite{Etnyre2006}.
\begin{lemma}[\cite{Etnyre2006}, ~Etnyre]\label{lem:etnyre}
    Suppose that $\S^3= \mathcal{A}ob(\Sigma,\Phi)$ and $(B,\pi)$ is an open book of $\S^3$ associated to $(\Sigma,\Phi)$. Let $\{\gamma_1,\gamma_2,\dots ,\gamma_n\}$ be a collection of simple closed curves on $\Sigma$ and let $\tau_{\gamma_i}$ denote a Dehn twist along $\gamma_i$. Then, the $3$-manifold 
\[
    M = \mathcal{A}ob(\Sigma, \displaystyle\prod_i\tau_{\gamma_i} \circ \Phi)
\]
    is obtained from $\S^3$ by surgery along a link $L=L_1\sqcup L_2 \sqcup \cdots \sqcup L_n$, where each $L_i$ is a copy of $\gamma_i$ such that they lie on distinct pages of the open book $(B,\pi)$. The surgery coefficient of each $L_i$ is $\pm 1$ with respect to the framing induced on $L_i$ by the page that contains $L_i$.
\end{lemma}
    
    In fact, below, we have another such lemma. First, recall that an embedded annulus in $\S^3$ is said to be \emph{standard} if the boundary components of the embedded annulus are unknots and unlinked. An embedded annulus is said to be a Hopf annulus if the boundary components are unknots with linking number $\pm 1$.
\begin{lemma}\label{lem:open_book_surgery}
    Suppose that $\S^3= \mathcal{A}ob(\Sigma,\Phi)$ with associated open book $(B,\pi)$ of $\mathbb{S}^3$. Let $\{L_1, L_2 ,\dots, L_n\}$ be a collection of disjoint simple closed curves on the page $\Sigma$ of the open book $(B,\pi)$ with the following properties
\begin{enumerate}
    \item each $L_i$ is an unknot in $\S^3$,
    \item for $1\le i \le k$, the tubular neighborhood of each $L_i$ in $\Sigma$ is a standard annulus in $\S^3$,
    \item for $k+1\le i \le n$, the tubular neighborhood of each $L_i$ in $\Sigma$ is a Hopf annulus in $\S^3$.
\end{enumerate}    
    Then the $3$-manifold $$M = \mathcal{A}ob(\Sigma, \tau_{L_1}^2\circ \cdots \circ \tau_{L_k}^2\circ \tau_{L_{k+1}} \circ \cdots \circ \tau_{L_n} \circ \Phi)$$ is obtained from $\S^3$ by surgery along the curves $L_i, 1\leq i\leq k$ with surgery coefficient $\pm \frac{1}{2}$ with resepect to Seifert framing and on the curves $L_i, k+1\leq i\leq n$ with surgery coefficient $\pm 1$ with respect to the framing induced by the page $\Sigma$ of the open book $(B,\pi)$. In fact, $M$ can be obtained from $\S^3$ by surgery on the link $$L=L_1\sqcup L_2\sqcup\cdots\sqcup L_n,$$ where $L_i$, $1\le i\le k$, are assigned surgery coefficients $\pm\frac{1}{2}$, while $L_i$, $k+1\le i\le n$, are assigned surgery coefficient $0$ or $\pm 2$, with respect to the Seifert framing.
\end{lemma}
\begin{figure}[ht]
    \centering
    \includegraphics[width=0.7\linewidth]{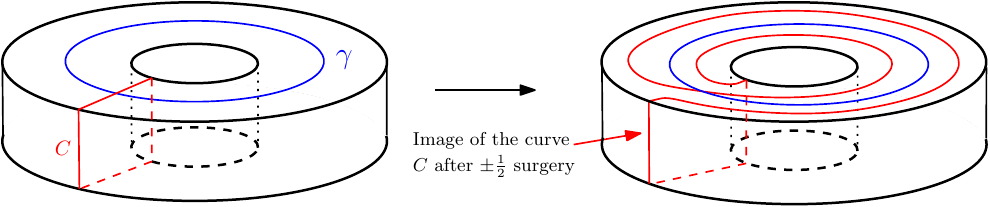}
    \caption{Effect of $\pm \frac{1}{2}$ surgery with respect to page framing.}
    \label{Open_book_after_surgery}
\end{figure}
\begin{remark}
    Note that, in Lemma~\ref{lem:open_book_surgery}, the monodromy presentation of the open book can be changed by reordering the Dehn twist factors in the composition $$\tau_{L_1}^2\circ \cdots \circ \tau_{L_k}^2\circ \tau_{L_{k+1}} \circ \cdots \circ \tau_{L_n},$$ since their supports are disjoint.

\end{remark}
    We now recall the notion of \emph{flexibility} from~\cite{HiroseYasuhara2008}.
\begin{definition}[$\phi$-Flexible Embedding]
    Let $f: \Sigma \hookrightarrow \mathbb{D}^4$ be a proper embedding of a compact, oriented surface $\Sigma$ with non-empty boundary in $\mathbb{D}^4$. If $f$ is  isotopic to an embedding $g: \Sigma \hookrightarrow \mathbb{D}^4$ such that 
    \begin{enumerate}
        \item $f(\Sigma) = g(\Sigma)$, and 
        \item $g^{-1} \circ f = \phi$,
    \end{enumerate} where $\phi : \Sigma \rightarrow \Sigma$ is a diffeomorphism of $\Sigma$, then we call the embedding $f$, a $\phi$-flexible embedding.
\end{definition}
    Equivalently, by the Isotopy Extension Theorem, a proper embedding $f\colon \Sigma\to \mathbb{D}^4$ is said to be $\phi$-flexible if there exists a diffeomorphism $\psi\colon \mathbb{D}^4 \to \mathbb{D}^4$ isotopic to the identity map such that $f^{-1}\circ \psi \circ f=\phi$. Now, we have the following straightforward lemma.
\begin{lemma}\label{lem:abs_open_book_embedding}
    Let $\phi$ be a diffeomorphism of a compact, oriented surface $\Sigma$ with non-empty boundary. If there exists a $\phi$-flexible embedding of $\Sigma$ in $\mathbb{D}^4$, then $M=\mathcal{A}ob(\Sigma,\phi)$ admits an abstract open book embedding into the trivial abstract open book of $\S^5=\mathcal{A}ob(\mathbb{D}^4,Id)$.
\end{lemma}
    The following lemma is well known and follows from the work of Hirose and Yasuhara~\cite{HiroseYasuhara2008},~\cite{Hirose2005}.

\begin{figure}[ht]
    \centering
    \includegraphics[width=0.8\linewidth]{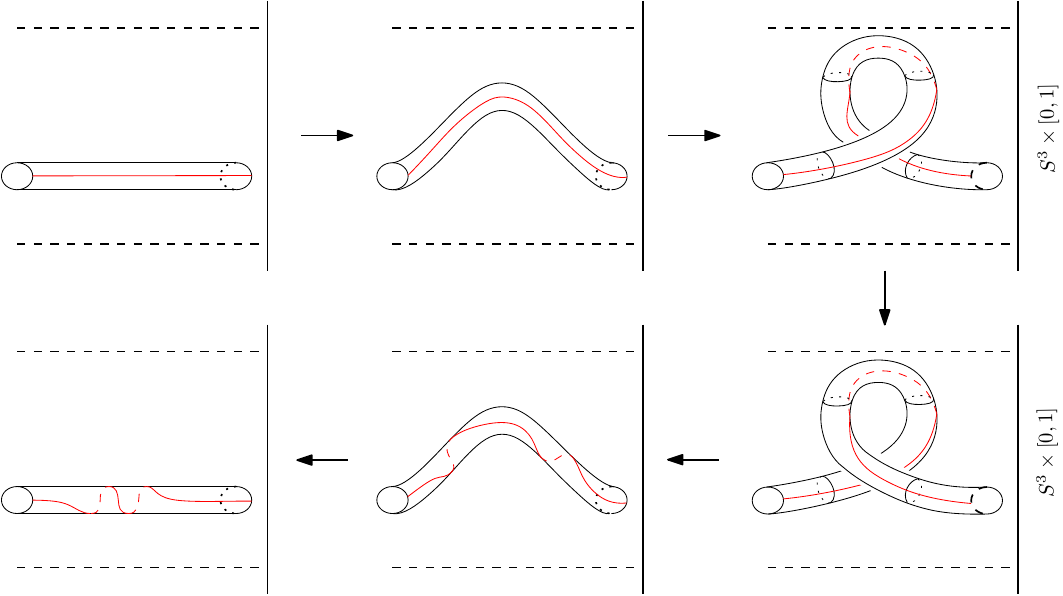}
    \caption{Tube trick.}
    \label{Tube trick}
\end{figure}

\begin{lemma}\label{lem:flex_curves}
    Let $\S^3\times [0,1]$ be a collar neighborhood of the boundary $\partial \mathbb{D}^4$ of $\mathbb{D}^4$ with $\S^3\times 0=\partial \mathbb{D}^4$. Let $f\colon \Sigma\to \mathbb{D}^4$ be a proper embedding of a compact, oriented surface $\Sigma$ with a non-empty boundary in $\mathbb{D}^4$ such that $f(\Sigma)\subset \S^3\times [0,1]$. Then, we have the following
\begin{enumerate}
    \item $\textbf{Hopf Annulus Trick:}$ Let $\gamma$ be a simple closed curve on $\Sigma$ such that $f(\gamma)\subset \S^3\times 1$ and a tubular neighborhood of $f(\gamma)$ in $f(\Sigma)$ is a Hopf annulus in $\S^3\times 1$, then $f$ is  $\tau_{\gamma}$-flexible.
    \item $\textbf{Tube Trick:}$ Let $\gamma$ be a simple closed curve on $\Sigma$ such that $f(\gamma)\subset \S^3\times 1$ and a tubular neighborhood of $f(\gamma)$ in $f(\Sigma)$ is a standard annulus in $\S^3\times 1$, then $f$ is $(\tau_{\gamma})^2$-flexible. 
\end{enumerate}
\end{lemma}
    For the pictorial description of the tube trick, see Figure~\ref{Tube trick}.
    The following lemma is straightforward, since every isotopy of $\mathbb{S}^3$ extends to an isotopy of $\mathbb{D}^4$.
\begin{lemma}\label{lem:phi-flexibility}
    Let $\S^3\times [0,1]$ be a collar neighborhood of the boundary $\partial \mathbb{D}^4$ of $\mathbb{D}^4$ with $\S^3\times 0=\partial \mathbb{D}^4$. Let $\Sigma$ be a fiber surface in $\S^3$ associated to an open book $(B,\pi)$. Let $(\Sigma,\phi)$ be the abstract open book associated to the open book  $(B,\pi)$ of $\S^3$. Consider an embedding of $\Sigma$ by inclusion in $\S^3\times 1$ and then by attaching the cylinders $\partial \Sigma \times [0,1]$ to the boundary of $\Sigma$ smoothly, we construct a smooth proper embedding $f$ of $\Sigma$ into $\mathbb{D}^4$. Then the embedding $f$ is $\phi$-flexible.
\end{lemma}
\bigskip

\section{Proof of Main Theorem}\label{sec:main}

    In this section, we discuss proof of Theorem \ref{thm:planar_page_main}.

\begin{proof}[\textbf{A Proof of Theorem~\ref{thm:planar_page_main}:}]
    Let $M$ be the given closed, oriented $3$-manifold. By Theorem~\ref{thm:binary_surgery_link}, we can obtain $M$ from $\S^3$ by surgery along a binary link $L=L_1\sqcup L_2\sqcup\cdots\sqcup L_n,$ where the surgery coefficient of $L_i$, for each $1\le i\le k$, is $\pm\frac{1}{2}$, and the surgery coefficient of $L_i$, for each $k+1\le i\le n$, is $0$. 

    By Theorem~\ref{thm:sinem}, there is a planar open book $(B^\prime,\pi^\prime)$ of $\S^3$ such that link $L$ lies on the fiber surface $\Sigma^1$ of $(B^\prime,\pi^\prime)$. We can appropriately perform a sequence of stabilizations of the open book $(B^\prime,\pi^\prime)$ to obtain an open book $(B,\pi)$ of $\S^3$ such that the fiber surface (page) $\Sigma$ associated to the book $(B,\pi)$ has the following properties:
\begin{enumerate}
    \item the fiber surface $\Sigma$ is a planar surface,
    \item  the link $L$ after isotoping in $\S^3$ lies on the fiber surface $\Sigma$, see Figure~\ref{Isotopy of the component of the link L after 1-handle attachment},
    \item the component of $L$ with surgery coefficient $\pm\frac{1}{2}$ has the framing with respect to page $\Sigma$ equal to the seifert framing,
    \item the component of $L$ with surgery coefficient $0$ has the framing with respect to page $\Sigma$ equal to $\pm1$ the Seifert framing.
\end{enumerate}

    To obtain the required fiber surface $\Sigma$, we first perform finitely many stabilizations of the fiber surface $\Sigma^1$ to obtain a fiber surface $\Sigma^\prime$ with the following property: For each component $L_i$ of $L$, there exists an arc $\alpha_i\subset\Sigma^\prime$ with one endpoint on $\partial\Sigma^\prime$ and the other on $L_i$, such that the interior of $\alpha_i$ is disjoint from the link $L$. 

    Suppose that, for some $j$, no such arc exists. This means that every arc joining $L_j$ to the boundary of $\Sigma^1$ must intersect at least one other component of $L$. Choose an arc $\alpha$ joining $L_j$ to a boundary component $B$ of $\Sigma^1$ such that the number of intersection points $|\operatorname{int}(\alpha)\cap L|$ is minimum and the interior of $\alpha$ is disjoint from $L_j$.

    We then stabilize the open book by attaching a $1$-handle to a single boundary component $B$ of the page $\Sigma^1$ near the point $\alpha\cap B$, ensuring that the stabilization curve $\gamma$ is disjoint from the link $L$ and that the attaching region does not contain the point $\alpha\cap B$, as shown in Figure~\ref{Isotopy of the component of the link L after 1-handle attachment (1)}. We denote the resulting stabilized fiber surface by $\Sigma^2$.

    Let $L_r$ be the first component of $L$ that we encounter while moving along $\alpha$ from $B$ towards $L_j$. We now isotope $L_r$ using a disk bounded by the stabilization curve $\gamma$, ensuring that the isotoped component remains on the stabilized fiber surface $\Sigma^2$, as shown in Figure~\ref{Isotopy of the component of the link L after 1-handle attachment (1)}. This isotopy removes the first intersection of $L_r$ with $\alpha$ without introducing any new intersections with $\alpha$.
\begin{figure}[ht]
    \centering
    \includegraphics[width=0.95\linewidth]{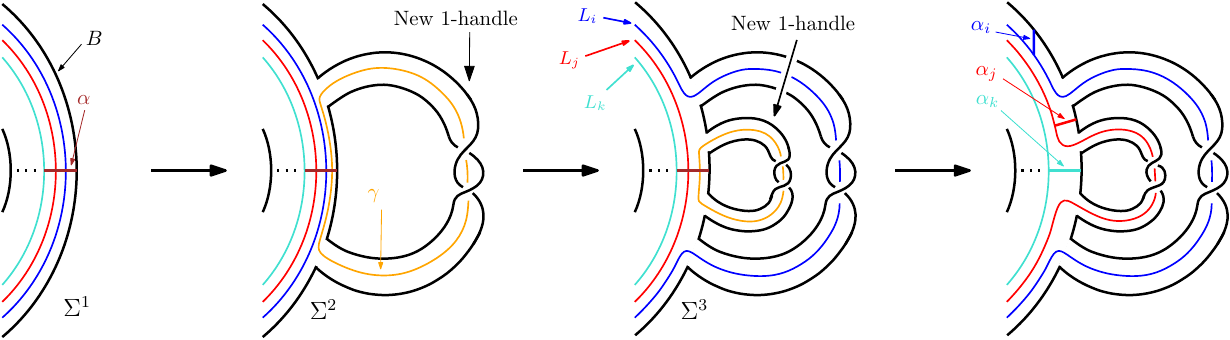}
    \caption{Local change in the component $L_i$ of the link L after isotopy.}
    \label{Isotopy of the component of the link L after 1-handle attachment (1)}
\end{figure}

    By repeating this procedure, we eventually obtain a stabilized fiber surface $\Sigma^l$, for some $l\in\mathbb{N}$, and an arc $\alpha_j$ joining $L_j$ to $\partial\Sigma^l$ such that the interior of $\alpha_j$ is disjoint from $L$.
    
    Since the link $L$ has finitely many components, we can repeat this construction for each component. Thus, after finitely many stabilizations, we obtain a fiber surface $\Sigma^{\prime}$ such that, for every component $L_j$ of $L$, there exists an arc $\alpha_j\subset\Sigma^{\prime}$ with one endpoint on $\partial\Sigma^{\prime}$ and the other on $L_j$, whose interior is disjoint from $L$.

    Now, to obtain a fiber surface $\Sigma$ with the desired framing, we again perform a sequence of stabilizations, each of which is obtained by attaching a $1$-handle to a single boundary component $B^\prime$ of the page $\Sigma^{\prime}$ parallel to a component $L_i$ of the link $L$, i.e., if there exists an arc $\alpha_i\subset\Sigma^{\prime}$ with one endpoint on $B^\prime$ and the other on $L_i$, whose interior is disjoint from $L$ as shown in Figure~\ref{Isotopy of the component of the link L after 1-handle attachment}. We choose the stabilization curve $\gamma$ such that it does not intersect the link $L$.  We then isotope $L_i$ using the disk bounded by the stabilization curve $\gamma$, as shown in Figure~\ref{Isotopy of the component of the link L after 1-handle attachment}, such that the isotoped component remains on the stabilized fiber surface. Note that each such a stabilization changes the page framing of the link component $L_i$ by $\pm 1$.
\begin{figure}[ht]
    \centering
    \includegraphics[width=0.8\linewidth]{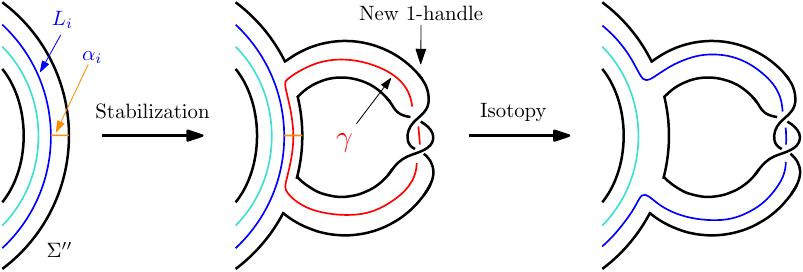}
    \caption{Local change in the component $L_i$ of the link L after isotopy.}
    \label{Isotopy of the component of the link L after 1-handle attachment}
\end{figure}

     Now $\Sigma$ is a fiber surface of the planar open book say $(B,\pi)$ of $\S^3$. Let $(\Sigma,\Phi)$ be the associated abstract open book to the open book $(B,\pi)$. We have the binary link $L\subset \Sigma$ with binary coefficients. If we perform surgery on the components of $L$ with assigned binary coefficients then the resulting manifold is the given manifold $M$. As in Lemma~\ref{lem:open_book_surgery}, the manifold $M$ admits an abstract open book $(\Sigma,\phi)$, where $$\phi=\tau_{L_1}^2\circ \cdots \circ \tau_{L_k}^2\circ \tau_{L_{k+1}} \circ \cdots \circ \tau_{L_n} \circ \Phi.$$
    
    Let $\S^3\times [0,1]$ be a collar neighborhood of the boundary $\partial \mathbb{D}^4$ of $\mathbb{D}^4$ with $\S^3\times 0=\partial \mathbb{D}^4$. Consider an embedding of $\Sigma$ by inclusion in $\S^3\times 1$  and then by attaching the cylinders $\partial \Sigma\times [0,1]$ smoothly we construct a smooth proper embedding $f$ of $\Sigma$ into $\mathbb{D}^4$. By Lemma~\ref{lem:flex_curves} the embedding $f$ is $\Phi$-flexible and by Lemma \ref{lem:phi-flexibility}, $f$ is $\tau_{L_i}^2$-flexible for $1\le i\le k$ as well as $\tau_{L_i}$-flexible for $k+1\le i\le n$. Hence, the embedding $f$ is $\phi$-flexible. Therefore, by Lemma~\ref{lem:abs_open_book_embedding}, $M=\mathcal{A}ob(\Sigma,\phi)$ admits an abstract open book embedding into $\S^5=\mathcal{A}ob(\mathbb{D}^4,Id)$. Thus the $3$-manifold admits a planar open book which open book embeds into the trivial open book of $\S^5$.
\end{proof}

\subsection{A Remark on Open Book Embedding}\mbox{}

   We have proved that every closed, orientable $3$-manifold $M$ admits a planar open book that open book embeds into the trivial open book of $\mathbb{S}^5$. This result relies on Theorem~\ref{thm:sinem} of Onaran. A similar result was established by Stallings~\cite{JRS}, which states that every link in $\mathbb{S}^3$ lies on a page of an open book of $\mathbb{S}^3$. More precisely, he proved the following theorem.

\begin{theorem}[\cite{JRS}, Stallings]\label{thm:stallings}
    Given a link $L$ in $\S^3$, there exists an unknot $K$ in $\S^3$ such that $L \sqcup K$ is the boundary of a fiber surface $\Sigma$ in $\S^3$. In particular, $L \sqcup K$ is a fibered link. Furthermore, we can choose $K$ such that the linking number of $K$ with any component of $L$ is a fixed pre-assigned integer.
\end{theorem}
    We would like to point out that using Stallings' result and the arguments given in the proof of Theorem~\ref{thm:planar_page_main}, we can establish a result analogous to Theorem~\ref{thm:planar_page_main}, although the resulting open book need not have a planar page. 
\medskip
\subsection{Embeddings of Open Books With Connected Binding}\label{subsec_stabi}\mbox{}

    Let $(B,\pi)$ be an open book decomposition of a closed, oriented $3$-manifold $M$ obtained in Theorem~\ref{thm:planar_page_main} with associated abstract open book $(\Sigma,\phi)$. Let $(\Sigma',\phi')$ be the stabilization of $(\Sigma,\phi)$. Then
\[
        (\Sigma',\phi')=(\Sigma\cup H^1,\phi\circ\tau_c),
\]
    where $H^1$ is a $2$-dimensional $1$-handle attached to $\Sigma$, and $\tau_c$ denotes the right- (or left-) handed Dehn twist along a simple closed curve $c \subset \Sigma^\prime$ that intersects the co-core of the $1$-handle $H^1$ exactly once.

    We now explain how to modify the embedding constructed in the proof of Theorem~\ref{thm:planar_page_main} to obtain a proper embedding of $\Sigma'$ into $\mathbb{D}^4$ which is $\phi^\prime$-flexible. Recall that, in the proof of Theorem~\ref{thm:planar_page_main}, we obtained a proper embedding of $\Sigma$ into $\mathbb{D}^4$ by first embedding $\Sigma$ into $\mathbb{D}^4$ via the standard inclusion and then smoothly attaching cylinders along the boundary components.

    To obtain an embedding of $\Sigma'$, we distinguish two cases depending on whether the $1$-handle $H^1$ is attached to two distinct boundary components of $\Sigma$ or to a single boundary component. In either case, we attach a smoothly embedded $1$-handle in $S^3\times\{1\}\subset\mathbb{D}^4$ to the corresponding boundary components of the embedded surface $\Sigma$ so that a regular neighborhood of the image of the curve $c$ becomes a Hopf annulus in $S^3\times\{1\}$. Finally, we smoothly attach the cylinders $\partial\Sigma'\times[0,1]$ along the boundary of $\Sigma'$ to obtain a smooth proper embedding $f:\Sigma'\hookrightarrow\mathbb{D}^4.$ Using arguments similar to those in the proof of Theorem~\ref{thm:planar_page_main}, we can prove that the embedding $f$ of $\Sigma'$ is $\phi'$-flexible.

    As an immediate consequence of the above construction, we obtain the Theorem~\ref{thm:connected_binding_main}, which states that every closed orientable $3$-manifold admits an open book with connected binding which open book embeds into the trivial open book of $\S^5$.

\section{Further discussions}\label{sec:fd}
    In the proof of Theorem~\ref{thm:binary_surgery_link}, by applying blow-ups, blow-downs, and Rolfsen twists, we can also obtain a surgery link in $\S^3$ for any closed, oriented $3$-manifold such that the surgery coefficients of its components belong to the set
\[
    \left\{0,\pm 2,\pm\frac{1}{2k}:k\in\mathbb Z\setminus\{0\}\right\}.
\]
    Using such surgery links, we can construct open book decompositions of any closed, oriented $3$-manifold $M$ that admit open book embeddings into the trivial open book of $\S^5$ as described in Theorems~\ref{thm:planar_page_main} and \ref{thm:connected_binding_main}.

    As we can see from Section \ref{sec:main}, there are infinitely many pairwise non-isotopic open book decompositions of $M$ that admit open book embeddings into the trivial open book of $\S^5$. This naturally leads to the following question: 
\begin{question}\label{que:all_open_book_embeddings_exists_or_not}
    Does every open book decomposition of a closed, orientable $3$-manifold embed into the trivial open book of $\mathbb{S}^5$?
\end{question}

    In~\cite{EtnyreOzbagci2008}, Etnyre and Ozbagci introduced three invariants of contact structures using open books: \emph{support genus}, \emph{binding number}, and \emph{support norm}, as described in Section~\ref{sec:intro}. 

    If the answer to Question~\ref{que:all_open_book_embeddings_exists_or_not} is not affirmative, then we can define a new invariant of closed, orientable $3$-manifolds using open books, which we call the \emph{embedding genus}.
\begin{definition}
    The \emph{\textbf{embedding genus}} of a closed, orientable $3$-manifold $M$ is the minimum genus among all open books with connected binding that admit an open book embedding into the trivial open book of $\mathbb{S}^5$.
\end{definition}
    We will denote it by $\mathcal{EG}(M)$.  We know that every closed, orientable $3$-manifold $M$ admits an open book decomposition with connected binding. Thus, we can define the \emph{minimal open book genus} of $M$ as follows.
\begin{definition}
    The \emph{\textbf{minimal open book genus}} of a closed, orientable $3$-manifold $M$ is the minimum genus of the page among all open book decompositions of $M$ with connected binding.
\end{definition}
    We denote the minimal open book genus of $M$ by $\mathcal{MG}(M)$. This naturally leads to the following question.
\begin{question}
    Is the embedding genus of $M$ equal to the minimal open book genus of $M$, i.e., $\mathcal{MG}(M)=\mathcal{EG}(M)$?
\end{question}
    For $\mathcal{EG}(M)\leq 2$, the equality $\mathcal{MG}(M)=\mathcal{EG}(M)$ holds. Since for $g\leq 1$, the surface $\Sigma_{g,1}$ admits a flexible proper embedding into $\mathbb{D}^4$, where $\Sigma_{g,n}$ denotes the compact orientable surface of genus $g$ and $n$ boundary components.
    
    By Theorem~\ref{thm:planar_page_main}, every closed, oriented $3$-manifold admits a planar open book decomposition that can be embedded into the trivial open book of $\mathbb{S}^5$. This motivates us to introduce another invariant of closed, orientable $3$-manifolds, which we call the \emph{embedding binding number}.
\begin{definition} 
    The \textbf{embedding binding number} of a closed, orientable $3$-manifold $M$ is defined to be the minimum number of boundary components of the page among all planar open book decompositions of $M$ that open book embeds into the trivial open book of $\mathbb{S}^5$.
\end{definition}
    We will denote it by $\mathcal{EBN}(M)$. The \emph{\textbf{minimal binding number}} of planar open books of $M$, denoted by $\mathcal{MBN}(M)$, is the minimum number of binding components of the page among all planar open book decompositions of $M$. This naturally leads to the following question.
\begin{question}
    Is the embedding binding number of $M$ equal to the minimal binding number of $M$, i.e., $\mathcal{EBN}(M)=\mathcal{MBN}(M)$?
\end{question}
    For $\mathcal{EBN}(M)\leq 2$, the equality $\mathcal{MBN}(M)=\mathcal{EBN}(M)$ holds. Since for $n\leq 2$, the surface $\Sigma_{0,n}$ admits a flexible proper embedding into $\mathbb{D}^4$.

\bigskip

    \noindent \emph{\textbf{Acknowledgements.}} Suhas Pandit is supported by MATRICS grant [MTR/2023 /000438] of Science and Engineering Research Board Government of India. Sanjay Kute is supported by NBHM Ph.D Scholarship [order number 0203/14/2021-R\&D-II/9396] India.

    \noindent \textbf{AI Declaration:} The free version of ChatGPT (OpenAI) was used solely for grammar checking and language improvement during the preparation of this manuscript. No AI tools were used for mathematical research, proofs, analysis, or the generation of scientific results.
\bibliographystyle{plain} 
\bibliography{references}
\end{document}